\documentclass[12pt,english,a4paper,final,reqno]{amsart}

\usepackage{showkeys}

\usepackage{geometry}
\usepackage[T1]{fontenc}
\usepackage[latin1]{inputenc}
\usepackage{babel}

\usepackage{setspace}
\usepackage{indentfirst}

\usepackage{hyperref}
\hypersetup{pdfpagemode=None,colorlinks=false,pdfpagelayout=OnePage,pdfstartview=FitV}

\makeatletter

\usepackage{verbatim}

\makeatother

\allowdisplaybreaks[1]
\usepackage{amssymb}
\usepackage{wasysym}
\usepackage{pifont}
\usepackage{stmaryrd}

\usepackage{graphicx}
\usepackage[all]{xy}
\def\dar[#1]{\ar@<2pt>[#1]\ar@<-2pt>[#1]}

\usepackage{amsthm}

\theoremstyle{plain}
\newtheorem{prop}{Proposition}[section]
\newtheorem{lem}[prop]{Lemma}

\newtheorem{cor}[prop]{Corollary}
\newtheorem{thm}[prop]{Theorem}

\newtheorem*{prop*}{Proposition}
\newtheorem*{lem*}{Lemma}
\newtheorem*{sublem*}{Sublemma}
\newtheorem*{cor*}{Corollaire}
\newtheorem*{thm*}{Th\'eor\`eme}
\newtheorem*{hypo*}{Hypothesis}
\newtheorem*{question*}{Question}
\newtheorem*{conjecture*}{Conjecture}
\newtheorem*{scholum*}{Scholum}
\newtheorem{defn}[prop]{Definition}
\newtheorem*{defn*}{D\'efinition}

\theoremstyle{definition}

\newtheorem{rmk}[prop]{Remark}

\newtheorem{ex}[prop]{Example}

\newtheorem*{con*}{Construction}
\newtheorem*{note*}{Note}
\newtheorem*{rmk*}{Remark}
\newtheorem*{rmks*}{Remarks}
\newtheorem*{ex*}{Example}
\newtheorem*{exs*}{Examples}

\theoremstyle{remark}

\newtheorem*{warning*}{Warning}
\newtheorem*{shortnote*}{Note}
\newtheorem*{claim*}{Claim}
\newtheorem*{axiom*}{Axiom}
\newtheorem*{example*}{Example}
\newtheorem*{examples*}{Examples}
\newtheorem*{remark*}{Remark}
\newtheorem*{remarks*}{Remarks}

\DeclareMathOperator{\idn}{id}

\DeclareMathOperator{\Aut}{Aut}
\DeclareMathOperator{\aut}{\mathfrak{aut}}

\DeclareMathOperator{\pr}{pr}

\DeclareMathOperator{\RE}{Re} 

\renewcommand{\epsilon}{\varepsilon}
\newcommand{\bphi}{\phi}
\newcommand{\vphi}{\varphi}

\newcommand{\beq}[1]{\begin{equation}\label{#1}}
\newcommand{\eeq}{\end{equation}}

\newcommand{\CC}{\mathbb{C}}
\newcommand{\RR}{\mathbb{R}}

\newcommand{\derlie}{\mathcal{L}} 
\newcommand{\liederivative}{L}

\newcommand{\ip}[2]{\langle #1,#2 \rangle} 
\newcommand{\lie}[2]{[#1,#2]} 
\newcommand{\schouten}[2]{[#1,#2]} 
\newcommand{\courant}[2]{\llbracket#1,#2\rrbracket} 
\newcommand{\dorfman}[2]{\courant{#1}{#2}+\DD\ip{#1}{#2}} 
\newcommand{\db}{\circ} 
\newcommand{\anchor}{\rho} 
\newcommand{\delbar}{\bar{\partial}}

\newcommand{\thalf}{\tfrac{1}{2}}
\newcommand{\half}{\frac{1}{2}}
\newcommand{\rond}{\circ} 
\newcommand{\cc}[1]{\overline{#1}} 

\newcommand{\gendex}[2]{\left\{ #1 \right\}_{#2}}

\newcommand{\cinf}{C^{\infty}} 
\newcommand{\sections}[1]{\Gamma(#1)}
\newcommand{\vf}{\mathfrak{X}} 
\newcommand{\df}{\Omega} 

\newcommand{\inclusion}{\hookrightarrow}

\newcommand{\isomorphism}{\simeq}

\newcommand{\diese}{^{\sharp}} 
\newcommand{\bemol}{^{\flat}}
\newcommand{\cmplx}{_{\mathbb{C}}} 
\newcommand{\inv}{^{-1}}

\newcommand{\DD}{\mathcal{D}}

\newcommand{\der}{\delta}
\newcommand{\ot}{\omega_t}
\newcommand{\otd}{\dot{\omega}_t}
\newcommand{\jtd}{\dot{J}_t}
\newcommand{\jt}{J_t}
\newcommand{\lt}{L_t}
\newcommand{\ltb}{\cc{L_t}}
\newcommand{\zt}{z_t}
\newcommand{\ztb}{\cc{z_t}}
\newcommand{\xt}{x_t}
\newcommand{\bt}{\beta_t}
\newcommand{\dlt}{d_{\lt}}
\newcommand{\ptd}{\dot{\pi}_t}
\newcommand{\Ot}{\Omega_t}
\newcommand{\Otd}{\dot{\Omega}_t}
\newcommand{\aoz}{A^{1,0}}
\newcommand{\aozt}{A^{1,0}_t}
\newcommand{\azo}{A^{0,1}}
\newcommand{\azot}{A^{0,1}_t}
\newcommand{\action}{\triangleright}
\newcommand{\proz}{\pr^{1,0}}
\newcommand{\cochains}[3]{C^{#1}(#2,#3)}
\newcommand{\dl}{d_L}
\newcommand{\yt}{Y_t}
\newcommand{\Ht}{H_t}
\newcommand{\Htd}{\dot{H}_t}
\newcommand{\eix}{e_i|_x}
\newcommand{\ejx}{e_j|_x}
\newcommand{\ekx}{e_k|_x}
\newcommand{\eitx}{e_i|_{tx}}
\newcommand{\ejtx}{e_j|_{tx}}

\newcommand{\TTM}{TM\oplus T^*M}

\newcommand{\liu}{MR1472888}

\newcommand{\tedcourant}{MR998124}
\newcommand{\gualtieri}{math.DG/0401221}
\newcommand{\hitchin}{MR2013140}
\newcommand{\barton}{math.DG/0603480}
\newcommand{\crainic}{math.DG/0412097}
\newcommand{\mackenzie}{MR1262213}
\newcommand{\roytenberg}{math.DG/9910078}
\newcommand{\newlander}{MR0088770}
\newcommand{\weinstein}{MR0286137}
\newcommand{\moser}{MR0182927}
\newcommand{\nest}{MR1913813}

\newcommand{\kodaira}{MR2109686}
\newcommand{\kohn}{MR0153030}
\newcommand{\cannas}{MR1747916}
\newcommand{\abraham}{MR515141}
\newcommand{\hp}{hp}
\newcommand{\severa}{letter}

\begin{document}
\title[]{Moser Lemma in Generalized Complex Geometry}
\author[]{Mathieu Sti\'enon}
\thanks{This work was supported by the European Union through the FP6 Marie Curie RTN ENIGMA (Contract number MRTN-CT-2004-5652) and by the E.S.I.~Vienna through a Junior Research Fellowship.}
\address{E.T.H.~Zürich, Departement Mathematik, 8092 Zürich, Switzerland}
\email{\href{mailto:stienon@math.ethz.ch}{stienon@math.ethz.ch}}
\begin{abstract}
We show how the classical Moser Lemma from symplectic geometry extends to generalized complex structures (GCS) on arbitrary Courant algebroids. For this, we extend the notion of Lie derivative to sections of the tensor bundle $(\otimes^i E)\otimes(\otimes^j E^*)$ with respect to sections of the Courant algebroid $E$ using the Dorfman bracket. We then give a cohomological interpretation of the existence of one-parameter families of GCS on $E$ and of flows of automorphims of $E$ identifying all GCS of such a family. 
In the particular cases of symplectic, we recover the results of Moser. Finally, we give a criterion to detect the local triviality of arbitrary GCS which generalizes the Darboux-Weinstein theorem.   
\end{abstract}


\maketitle


\section{Introduction}

The classical Moser lemma  for symplectic manifolds \cite{\moser}
describes a cohomological condition for two symplectic
structures to be equivalent. It can be stated as follows.
For two symplectic forms $\Omega_0$ and  $\Omega_1$ on a 
compact symplectic manifold $M$, if there exists a smooth 
one-parameter family $\omega_t$ of symplectic forms on $M$, 
all with the same periods and such that $\Omega_0=\omega_0$ 
and $\Omega_1=\omega_1$, then there is a global diffeomorphism 
$\varphi$ of $M$, diffeotopic to the identity and such that 
$\Omega_1=\varphi^*\Omega_0$.

Recently there has been increasing interest in 
generalized complex structures \cite{\hitchin, \gualtieri}, which
comprise both symplectic and complex structures as
special cases. It is natural to ask if Moser's lemma
extends to generalized complex geometry. The aim of
the present paper is to give an affirmative answer
to this question.
While up to now, most of the research on generalized
complex structures focused on exact Courant algebroids
$TM\oplus T^*M$ \cite{\severa}, in this paper,
we work on generic Courant algebroids.
A generalized complex structure on a Courant algebroid 
$E$ is a bundle map $J: E\to E$ satisfying $J^2=-1$, which is orthogonal 
with respect to the symmetric pairing and whose 
Nijenhuis torsion vanishes \cite{\barton}.

As is well known, the proof of the classical
Moser lemma involves the Lie derivative of
symplectic forms with respect to a vector field.
Here, as first step, we introduce a concept of 
Lie derivative of generalized complex strutures with respect to sections of a Courant algebroid. Indeed, we introduce the Lie derivative 
of sections of the tensor bundle 
$(\otimes^i E) \otimes (\otimes^j E^*)$ with respect to sections of the 
Courant algebroid $E$. We hope this construction will 
be of independent interest in the future. 
Such a Lie derivative can be defined exactly as in the classical
case. Namely, via the Dorfman bracket, one can think of $\sections{E}$ as
a subset of $\aut(E)$, the infinitesimal
automorphisms of the Courant algebroid $E$.
Thus any section of $E$ generates a local flow of automorphims of the 
Courant algebroid $E$. We define
the Lie derivative of any section of $(\otimes^i E)\otimes(\otimes^j E^*)$
as the time-derivative at $t=0$ of the pull back of this section
by the flow. For any $X,Y\in\sections{E}$, the Lie derivative 
$\derlie_X Y$ is simply the Dorfman bracket $X\db Y$. The Lie derivative of a 
bundle map $J: E\to E$, seen as a section of $E^*\otimes E$, 
is given by $(\derlie_XJ)(Y)=\derlie_X J(Y)-J(\derlie_X Y)$, $\forall X, Y\in \sections{E}$. This formula reduces exactly to the usual Lie derivative of
the symplectic form when $J$ corresponds to a symplectic structure
and $X$ is a vector field.

With the help of the Lie derivative, we are able to translate 
the infinitesimal isomorphism condition for $J_t$: 
\[ \jtd+\derlie_{\xt}\jt=0, \quad \forall t \text{ with } \xt\in\sections{E} \]
as the exactness, at every time $t$, of a family of $t\to\otd$ of Lie algebroid 2-cocycles relative to the one-parameter family of Lie algebroids 
$t\mapsto L_t$, where $L_t$ is the $+i$-eigenbundle of
$J_t$. Hence we recover exactly the same situation as in
the classical Moser lemma context \cite{\moser}.

As first examples, we consider symplectic
and complex structures on a Lie alebroid $A$.
The Courant algebroid $E$ involved is the double $A\oplus A^*$ 
of $A$. When $A$ is the tangent bundle Lie algebroid $TM$,
we recover the usual Moser lemma for symplectic manifolds and
the classical result of Kodaira \cite{\kodaira} respectively.
Other examples of generalized complex strutures
are given by Hamitonian operators in the sense of Liu-Weinstein-Xu
\cite{\liu}. We describe conditions when such operators induce isomorphic generalized complex structures.
Holomorphic Poisson structures are a special case.

As another application, we give a Darboux-Weinstein style theorem for
generalized complex strutures. More precisely, 
we describe local cohomological conditions which guarantee
the local triviality of a generalized
complex structure on the standard Courant algebroid. 
This cohomological condition is always satisfied for symplectic manifolds.
On the other hand, for an integrable complex structure, this condition
should be related to Kohn's proof \cite{\kohn} of the Newlander-Nirenberg theorem \cite{\newlander}.

\subsection*{Acknowledgements} The author is grateful to Alan Weinstein for the encouraging discussion they had in Berkeley at the early stage of this project and to Alberto Cattaneo, Vasiliy Dolgushev, Camille Laurent-Gengoux, Giovanni Felder, Pierre Schapira,  Boris Tsygan, Aissa Wade and Ping Xu for other ones; also to the E.S.I.~Vienna and I.H.P.~Paris for their hospitality while parts of this paper were written.

\section{Automorphisms of Courant algebroids}

\begin{defn}[\cite{\liu}]
A Courant algebroid is a triple consisting of
a vector bundle $E\to M$ equipped with 
a non degenerate symmetric bilinear form $\ip{\cdot}{\cdot}$,
 a skew-symmetric bracket $\courant{\cdot}{\cdot}$ on $\sections{E}$,
 and a smooth bundle map $E\xrightarrow{\anchor}TM$ called the anchor.
These induce a natural differential operator $\DD :\cinf(M)\to\sections{E}$
defined by \[ \ip{\DD f}{a}=\half \anchor(a)f \] for all $f \in \cinf(M)$ and $a\in\sections{E}$.

These structures must be compatible in the following sense:
 $\forall a,b,c\in\sections{E}$ and $\forall f,g\in \cinf(M)$,
\begin{align} 
& \anchor(\courant{a}{b})=\lie{\anchor(a)}{\anchor(b)}, \notag \\ 
& \courant{\courant{a}{b}}{c}+\courant{\courant{b}{c}}{a}
+\courant{\courant{c}{a}}{b}=\tfrac{1}{3}\DD\big(\ip{\courant{a}{b}}{c}+\ip{\courant{b}{c}}{a}+\ip{\courant{c}{a}}{b}\big), \notag \\ 
& \courant{a}{f b}=f\courant{a}{b}+\big(\anchor(a)f\big)b-\ip{a}{b}\DD f, 
\notag \\ 
& \anchor\rond\DD=0,\text{ i.e. }\ip{\DD f}{\DD g}=0, \notag \\ 
& \anchor(a)\ip{b}{c}=\ip{\courant{a}{b}+\DD 
\ip{a}{b}}{c}+\ip{b}{\courant{a}{c}+\DD \ip{a}{c}}. \label{c5} 
\end{align}
\end{defn}

\begin{ex}[\cite{\tedcourant}]
Given a smooth manifold $M$, the bundle $\TTM\to M$ carries a natural Courant algebroid structure,
where the anchor is the projection onto the tangent component and the pairing
and bracket are given, respectively, by
\begin{gather*}
\ip{X+\xi}{Y+\eta}=\tfrac{1}{2}\big(\xi(Y)+\eta(X)\big), \\
\courant{X+\xi}{Y+\eta}=\lie{X}{Y}+\liederivative_{X}\eta-\liederivative_{Y}\xi+\tfrac{1}{2}d\big(\xi(Y)-\eta(X)\big), 
\end{gather*}
$\forall X,Y \in \vf(M)$, $\forall \xi,\eta \in \df^1(M)$.
\end{ex}

The non-symmetric law  
\[ a\db b:=\courant{a}{b}+\DD \ip{a}{b} \]
is called \emph{Dorfman bracket}
\cite{\roytenberg}. The following result is due 
to Roytenberg \cite{\roytenberg}.

\begin{prop}
\label{pro:2.2}
For all $a,b,c\in\sections{E}$ and $f\in\cinf(M)$, one has: 
\begin{align*}
& \DD f\db a=0, &
& a\db \courant{b}{c}=\courant{a\db b}{c}+\courant{b}{a\db c}, \\
& a\db (fb)=f(a\db b)+ (\rho(a)f)b, & 
& \rho (a) \ip{b}{c}=\ip{a\db b}{c}+\ip{b}{a\db c}. 
\end{align*}
\end{prop}

Any section of the dual bundle $E^*\to M$ can be seen 
as a fiberwise linear function on $E$ and vice versa. 
In other words, $\sections{E^*}$ is naturally a subspace of $\cinf(E)$. 
Recall that an infinitesimal automorphism of
a vector bundle $E\stackrel{\pi}{\to} M$ corresponds
exactly to a vector field on $E$ (i.e. a derivation of $\cinf(E)$) 
under which the subspaces $\pi^*\cinf(M)$ and $\sections{E^*}$ 
are stable \cite[Proposition 2.2]{\mackenzie}.
The latter is equivalent to a covariant differential
operator on $E^*$ \cite{\mackenzie}, i.e. a pair of differential
operators $\der^0:\cinf(M)\to\cinf(M)$ and
$\der^1:\sections{E^*}\to\sections{E^*}$ satisfying 
\begin{align*}
\der^0(fg) &=f\der^0(g)+\der^0(f)g, \qquad \forall f,g\in\cinf(M) \\
\der^1(f\alpha) &=f\der^1(\alpha)+\der^0(f)\alpha, \qquad \forall \alpha\in\sections{E^*},\;\forall f\in\cinf(M) 
.\end{align*}

From now on, we assume that $E$ is a Courant algebroid.
By $\aut(E)$, we denote the Lie algebra of infinitesimal 
automorphism of the Courant algebroid $E$. The symmetric pairing 
$\ip{\cdot}{\cdot}$ identifies $\sections{E}$ with $\sections{E^*}$ 
by $x\mapsto\ip{x}{\cdot}$.
Therefore, one easily recovers the following theorem of 
Roytenberg \cite{\roytenberg}, which he proved using
super-geometry.

\begin{prop} 
\label{pro:2.3}
The  Lie algebra  $\aut(E)$ consists of those covariant differential
operators $\der=(\der^0, \der^1)$ on $E$, where 
 $\der^0: \cinf(M)\to \cinf(M)$ and
$\der^1: \sections{E}\to \sections{E}$, satisfying the additional
properties:
\[ \der^0\ip{x}{y}=\ip{\der^1 x}{y}+\ip{x}{\der^1 y} \] and \[ \der^1\courant{x}{y}=\courant{\der^1 x}{y}+\courant{x}{\der^1 y} ,\]
for all $x,y\in\sections{E}$.
\end{prop}

\section{Lie derivatives}
\label{liederiv}

For any $z\in \Gamma (E)$, define 
\[ \der^0_z: \cinf(M)\to \cinf(M) \qquad \text{and} \qquad 
\der_z^1: \sections{E}\to \sections{E} \] by 
\[ \der^0_z(f)=\anchor(z)(f) \qquad\text{and}\qquad \der_z^1(x)=z\db x, \qquad\forall
f\in\cinf(M),\;\forall x\in\sections{E} .\]
It follows from Propositions \ref{pro:2.2} and \ref{pro:2.3} 
that $\der_z=(\der^0_z, \der_z^1)$ is an infinitesimal
automorphism of the Courant algebroid $E$, i.e. $\der_z\in\aut(E)$. 

By $\bphi_t$ we denote the (local) flow 
\[ \xymatrix{E \ar[r]^{\bphi_t} \ar[d]_{\pi} & E \ar[d]^{\pi} \\ 
M \ar[r]_{\vphi_t} & M} \]
generated by the vector field on $E$ corresponding to $\der_z$. 
The 1-parameter group $t\mapsto\bphi_t$ acts on 
the space of sections of $E\to M$ by: 
\[ \bphi_t^*:\sections{E}\to\sections{E}:\sigma\mapsto\bphi_t\inv\rond\sigma\rond\vphi_t .\]

By abuse of notation, we use the same symbol
$\bphi_t$ (resp. $\bphi_t^*$) to denote the induced flow on the tensor
bundles $E_j^i=(\otimes^i E)\otimes(\otimes^j E^*)$ ($i,j\in\gendex{0,1,2,\dots}{}$) 
(resp. the induced action on the spaces of sections of the $E_j^i$'s). 

For any section $\sigma\in\sections{E_j^i}$, define
$\derlie_z\sigma\in\sections{E_j^i}$ by
\[ \derlie_z\sigma=\left.\tfrac{d}{d\tau}\bphi_{\tau}^*\sigma\right|_{\tau=0} \]
(see \cite[Theorem~2.2.20]{\abraham}).
Thus we have the usual identity:
\[ \left.\tfrac{d}{d\tau}\bphi_{\tau}^*\sigma\right|_{\tau=t}=\bphi_t^*(\derlie_z\sigma) .\]

In the following proposition, we give a list of important
properties of this Lie derivative, which will be useful
in the future discussion.

\begin{prop}
For all $f,g\in\cinf(M)$ and $x,y,z\in\sections{E}$, we have: 
\begin{align*}
& \derlie_z f=\rho(z) f, &
& \derlie_z\ip{x}{y}=\ip{\derlie_z x}{y}+\ip{x}{\derlie_z y}, \\
& \derlie_z x=z\db x, &
& \derlie_z\courant{x}{y}=\courant{\derlie_z x}{y}+\courant{x}{\derlie_z y}, \\
& \derlie_{\DD f}x=0, &
& \derlie_{fx} y=f\;\derlie_x y-(\derlie_y f)\;x+2\ip{x}{y}\DD f, \\  
& \derlie_{\courant{x}{y}}=\lie{\derlie_x}{\derlie_y}. && 
\end{align*}
Moreover, 
\begin{align*}
& \derlie_z(\sigma\otimes\tau)=\derlie_z\sigma\otimes\tau
+\sigma\otimes\derlie_z\tau \text{ for all } \sigma,\tau\in\oplus_{i,j}E^i_j, \\
& \text{and } \lie{\der}{\derlie_z}=\derlie_{\der^1 z} \text{ for all } \der\in\aut(E). 
\end{align*}
\end{prop}

In particular, if $J: E\to E$ is a bundle map over the identity $M\xrightarrow{\idn}M$, i.e. $J\in\sections{E^*\otimes E}$,
then \[ (\derlie_{x}J) (y)=\derlie_{x} (J(y))-J(\derlie_{x} y), \qquad\forall x,y\in\sections{E} .\]

\begin{ex}
For the standard Courant algebroid structure on $\TTM$, one has 
\[ \derlie_{X+\xi}(Y+\eta)=L_X (Y+\eta)-i_Y d\xi, 
\qquad \forall X,Y\in\vf(M),\;\forall \xi,\eta\in\df^1(M), \]
where $L$ denotes the usual Lie derivative.
\end{ex}

\begin{prop}
Let $E$ denote the standard Courant algebroid structure on $\TTM$. 
Then, for all $X\in\vf(M)$ and $\xi\in\df^1(M)$, the
1-parameter groups integrating the infinitesimal automorphisms 
$\der_X,\der_\xi,\der_{X+\xi}\in\aut(E)$ are given respectively by
\begin{gather*}
 e^{t\delta_X}(Y+\eta)= (\vphi_t)_*Y+(\vphi\inv_t)^*\eta, \\
 e^{t\delta_\xi}(Y+\eta)= Y+\eta+t\;i_Y\,d\xi, \\
 e^{t\delta_{X+\xi}}(Y+\eta)= (\vphi_t)_*Y+(\vphi\inv_t)^*\eta+\int_0^t(\vphi\inv_{t-\tau})^*\big(i_{\vphi_{\tau*}Y}\;d\xi\big)\;d\tau,
\end{gather*}
where $t\mapsto\vphi_t$ denotes the flow of the vector field $X$, 
i.e. satisfies 
$$ X_{\vphi_t(m)}=\left.\tfrac{d}{dt}\vphi_{\tau}(m)\right|_{\tau=t} .$$
In case the section $X_t+\xi_t\in\sections{E}$ is time-dependent, its flow
\[ \xymatrix{
\TTM \ar[d] \ar[r]^{\bphi_{a,b}} & \TTM \ar[d] \\
M \ar[r]_{\vphi_{a,b}} & M
} \] from time $a$ to $b$ is given by
\beq{eq:flow}
\bphi_{a,b}(Y+\eta)= (\vphi_{a,b})_*Y+(\vphi_{b,a})^*\eta+\int_a^b(\vphi_{b,\tau})^*\big(i_{(\vphi_{a,\tau})_*Y}\;d\xi_{\tau}\big)\;d\tau
,\eeq
where $\vphi_{a,b}$ denotes the flow of the time-dependent 
vector field $X_t$ from time $a$ to $b$.
\end{prop}

\begin{proof}
From \eqref{eq:flow} it follows that 
\[ \bphi_{b,c}\rond\bphi_{a,b}=\bphi_{a,c} \qquad\text{and}\qquad \bphi_{a,a}=\idn_E .\]
Since 
\[ \left.\tfrac{d}{dt}\bphi\inv_{a,a+\epsilon}\rond(Y+\eta)\rond\vphi_{a,a+\epsilon}\right|_{\epsilon=0}
=\left.\big(L_{X_t}(Y+\eta)-i_Y\;d\xi_t\big)\right|_{t=a} ,\]
the flow of $\der_{X_t+\xi_t}$ is indeed given by \eqref{eq:flow}.
\end{proof}

\section{Main theorem}

Let $E$ be a Courant algebroid on a smooth manifold $M$.
And let 
\[ \xymatrix{
E \ar[d] \ar[r]^J & E \ar[d] \\ 
M \ar[r]_{\idn} & M 
} \]
be a vector bundle map such that $J^2=-\idn$. 
Then the complexification $E\cmplx:=E\otimes\CC$ 
--- with the extended $\CC$-linear Courant algebroid structure --- 
decomposes as the direct sum $L\oplus\cc{L}$
of the eigenbundles of $J$. 
Here $L$ is associated to the eigenvalue $+i$ and its complex conjugate $\cc{L}$ to $-i$.
The bundle map $J$ is called a \emph{generalized complex structure} if $J$ is orthogonal with respect to $\ip{\cdot}{\cdot}$ --- this forces $L$ and $\cc{L}$ to be isotropic --- and the spaces of sections $\sections{L}$ and $\sections{\cc{L}}$ are closed under the Courant bracket, or equivalently, 
$J$ is ``integrable'': 
\[ \courant{Jx}{Jy}-\courant{x}{y}-J\big(\courant{Jx}{y}+\courant{x}{Jy}\big)=0, \qquad \forall x,y\in\sections{E} .\]

We refer the reader to \cite{\barton} for more details
and to \cite{\gualtieri} for the case of generalized complex structures
on the standard Courant algebroid $\TTM$.

The group of automorphisms of the Courant algebroid acts on its set of generalized complex structures by: 
\[ \phi_* J=\phi\inv\rond J\rond\phi ,\] where $\phi\in\Aut(E)$.
Two generalized complex structures $J_0$ and $J_1$
are said to be \emph{isomorphic} if there exists a
Courant algebroid automorphism $\phi$ such that $\phi_* J_1=J_0$.

Given a smooth family $t\mapsto\jt$ of generalized complex structures on $E$, 
we define a 2-form $\ot\in\sections{\wedge^2 E^*}$ by
\[ \ot(x,y)=\ip{x}{\jt y}, \qquad \forall x,y\in\sections{E} .\]
It is clear that $\ot$ is skew-symmetric.
By $L_t$ we denote the $+i$ eigenbundle 
of $J_t: E\cmplx \to E\cmplx$. 
Thus $L_t$ is a Dirac structure on
$E\cmplx$. Indeed $(L_t,\cc{L_t})$ is a (complex) Lie
bialgebroid \cite{\mackenzie}.

The restriction of the time-derivative of $\ot\in\sections{\wedge^2 E^*}$ 
to the subbundle $L_t$ defines a 2-form $\otd\in\sections{\wedge^2\lt^*}$:
\[ \otd(v,w)=\ip{v}{\jtd w},\qquad\forall v,w\in\sections{\lt} .\]

\begin{rmk}
Differentiating $\jt \jt=-\idn$, we obtain $\jtd\jt+\jt\jtd=0$. Therefore, the map $\jtd$ swaps the eigenbundles $\lt$ and $\ltb$.
\end{rmk}

Recall that there exists a canonical differential complex 
\[ \cdots\xrightarrow{d_{\lt}}\sections{\wedge^{\bullet-1}\lt^*}\xrightarrow{d_{\lt}}\sections{\wedge^\bullet \lt^*}
\xrightarrow{d_{\lt}}\sections{\wedge^{\bullet+1} \lt^*}\xrightarrow{d_{\lt}}\cdots \]
associated to the Lie algebroid $\lt$ \cite{\cannas}.
The coboundary operator is given by 
\begin{multline*}
\big(d_{\lt}\alpha\big)(x_0,x_1,\cdots,x_n)=\sum_{i=0}^n (-1)^i \big(\anchor x_i\big)
 \alpha(x_0,\cdots,\widehat{x_i},\cdots,x_n) \\
+ \sum_{i<j} (-1)^{i+j} \alpha(\courant{x_i}{x_j},x_0,\cdots,\widehat{x_i},\cdots,\widehat{x_j},\cdots,x_n)
, \end{multline*}
where $\alpha\in\sections{\wedge^n \lt^*}$ and $x_0,\dots,x_n\in\sections{\lt}$.
The cohomology of this complex is the \emph{Lie algebroid cohomology} of $\lt$.

\begin{prop} $d_{\lt}\otd=0$, $\forall t$ \end{prop}

\begin{proof}
\fbox{Step 1} Since $\jt$ is integrable $\forall t$, we have 
\[ \courant{\jt v}{\jt w}-\courant{v}{w}-\jt\big(\courant{\jt v}{w}+\courant{v}{\jt w}\big)=0 .\]
Differentiating w.r.t. $t$, we get 
\[ \courant{\jtd v}{\jt w}+\courant{\jt v}{\jtd w}-\jtd\big(\courant{\jt v}{w}+\courant{v}{\jt w}\big)
-\jt\big(\courant{\jtd v}{w}+\courant{v}{\jtd w}\big)=0 .\]
Now, taking $v,w\in\sections{\lt}$, the above relation becomes 
\[ i\;\courant{\jtd v}{w}+i\;\courant{v}{\jtd w}-2i\;\jtd\courant{v}{w}
-\jt\big(\courant{\jtd v}{w}+\courant{v}{\jtd w}\big)=0 \]
and, multiplying by $(-i)$, we obtain 
\[ \tfrac{1+i\jt}{2} \big(\courant{\jtd v}{w}+\courant{v}{\jtd w}\big) = \jtd\courant{v}{w} .\]
Hence, since $\lt$ is isotropic, 
\begin{equation} 
\ip{z}{\courant{\jtd v}{w}+\courant{v}{\jtd w}-\jtd\courant{v}{w}}=0,\qquad\forall v,w,z\in\sections{\lt} 
.\label{eq:21}\end{equation}
\fbox{Step 2} By definition, 
\[ \begin{split} 
\big(d_{\lt}\otd\big)(v,w,z) =& \anchor(v)\otd(w,z)-\anchor(w)\otd(v,z)+\anchor(z)\otd(v,w) \\ 
& -\otd(\courant{v}{w},z)+\otd(\courant{v}{z},w)-\otd(\courant{w}{z},v) \\
=& \anchor(v)\ip{w}{\jtd z}-\anchor(w)\ip{v}{\jtd z}+\anchor(z)\ip{v}{\jtd w} \\ 
& -\ip{\courant{v}{w}}{\jtd z}+\ip{\courant{v}{z}}{\jtd w}-\ip{\courant{w}{z}}{\jtd v}
.\end{split} \]
Now unfold the second and third terms according to \eqref{c5}:
\[ \begin{split} 
\big(d_{\lt}\otd\big)(v,w,z) 
=&\anchor(v)\ip{w}{\jtd z}-\ip{\dorfman{w}{v}}{\jtd z}-\ip{v}{\dorfman{w}{\jtd z}} \\
&+\ip{\dorfman{z}{v}}{\jtd w}+\ip{v}{\courant{z}{\jtd w}}+\ip{v}{\dorfman{z}{\jtd w}} \\
&-\ip{\courant{v}{w}}{\jtd z}+\ip{\courant{v}{z}}{\jtd w}-\ip{\courant{w}{z}}{\jtd v}
;\end{split} \]
simplify keeping in mind that $\lt$ is isotropic:
\[ \begin{split} 
\big(d_{\lt}\otd\big)(v,w,z)=&\anchor(v)\ip{w}{\jtd z}-\ip{v}{\courant{w}{\jtd z}}-\ip{v}{\DD\ip{w}{\jtd z}}
\\ &+\ip{v}{\courant{z}{\jtd w}}+\ip{v}{\DD\ip{z}{\jtd w}}-\ip{\courant{w}{z}}{\jtd v}
;\end{split} \]
and use the definition of $\DD$ to get 
\[ \begin{split} 
\big(d_{\lt}\otd\big)(v,w,z) =& 
\ip{v}{\courant{\jtd z}{w}}+\ip{v}{\courant{z}{\jtd w}}
+\ip{\jtd v}{\courant{z}{w}} \\ & +\anchor(v)\ip{w}{\jtd z}
-\thalf\anchor(v)\ip{w}{\jtd z}+\thalf\anchor(v)\ip{\jtd w}{z} \\
=& \ip{v}{\courant{\jtd z}{w}}+\ip{v}{\courant{z}{\jtd w}}\ip{\jtd v}{\courant{z}{w}} \\ 
& +\thalf\anchor(v)\big(\ip{\jtd w}{z}+\ip{w}{\jtd z}\big) 
.\end{split} \]
But, differentiating $\ip{\jt x}{\jt y}=\ip{x}{y}$ w.r.t. $t$, we get 
\[ \ip{\jtd x}{\jt y}+\ip{\jt x}{\jtd y}=0 .\]
Hence, if $x,y\in\sections{\lt}$, one has 
\[ \ip{\jtd x}{y}+\ip{x}{\jtd y}=0 .\]
Therefore, 
\[ \big(d_{\lt}\otd\big)(v,w,z)=\ip{v}{\courant{\jtd z}{w}+\courant{z}{\jtd w}-\jtd\courant{z}{w}} \] 
and the result follows from \eqref{eq:21}.
\end{proof}

\begin{lem} \label{lem:bt}
Let $t\mapsto\jt$ ($t\in[0,1]$) be a smooth 
family of generalized complex structures 
on a Courant algebroid $E$.
There exists a smooth path $t\mapsto\bt$ in $\sections{\lt^*}$ such that 
\beq{eq:wt} \otd=d_{\lt}\bt, \qquad \forall t \eeq
 if, and only if,
\[ \jtd+\derlie_{\xt}\jt=0, \qquad \forall t \]
where $\xt=\tfrac{1}{2i}(\zt-\ztb)\in\sections{E}$ is
 the imaginary part of the unique section
$\zt\in\sections{\lt}$ such that $\bt=\ip{\ztb}{\cdot}$.
\end{lem}

\begin{proof}
For $v,w\in\sections{\lt}$, one has $\otd(v,w)=\ip{v}{\jtd w}$ and 
\[ \begin{split}
\big(\dlt\bt\big)(v,w)&=\anchor(v)\bt(w)-\anchor(w)\bt(v)-\bt(\courant{v}{w}) \\ 
&=\anchor(v)\ip{\ztb}{w}-\anchor(w)\ip{\ztb}{v}-\ip{\ztb}{\courant{v}{w}} \\
&=2\ip{\DD\ip{\ztb}{w}}{v}-(\ip{\derlie_w\ztb}{v}+\ip{\ztb}{\derlie_w v})-\ip{\ztb}{\courant{v}{w}} \\
&=2\ip{\DD\ip{\ztb}{w}}{v}-\ip{\dorfman{w}{\ztb}}{v}-\ip{\ztb}{\courant{w}{v}}-\ip{\ztb}{\courant{v}{w}} \\ 
&=\ip{v}{\dorfman{\ztb}{w}}=\ip{v}{\derlie_{\ztb}w} 
.\end{split} \]
Hence 
\[ \begin{split} 
& \otd=\dlt\bt \\
\Leftrightarrow \qquad & \ip{v}{\jtd w-\derlie_{\ztb}w}=0,\quad\forall v\in\sections{\lt} \\ 
\Leftrightarrow \qquad & \jtd w-\derlie_{\ztb}w\in\sections{\lt} \\ 
\Leftrightarrow \qquad & \big(\jtd w\big)_{\ltb}=\big(\derlie_{\ztb}w\big)_{\ltb} 
,\end{split} \] for all $w\in\sections{\lt}$.
But $\jtd w\in\sections{\cc{\lt}}$ as $w\in\sections{\lt}$. 
Therefore, \begin{equation} \jtd w=\big(\derlie_{\ztb}w\big)_{\ltb}, \qquad\forall w\in\sections{\lt}  
.\label{eq:23}\end{equation}
Since $\sections{\lt}$ is $\courant{\cdot}{\cdot}$-closed and $\lt$ is $\ip{\cdot}{\cdot}$-isotropic, 
\[\derlie_{\zt}w=\dorfman{\zt}{w}=\courant{\zt}{w}\in\sections{\lt} \]
and \begin{equation} \big(\derlie_{\zt}w\big)_{\ltb}=0 ,\label{eq:24}\end{equation} for all $w\in\sections{\lt}$.
Substracting \eqref{eq:24} from \eqref{eq:23}, we get 
\begin{multline*} 
\jtd w=\big(\derlie_{(\ztb-\zt)}w\big)_{\ltb}=-2i\;\big(\derlie_{\xt}w\big)_{\ltb} 
=-2i\;\tfrac{1+i\jt}{2}\big(\derlie_{\xt}w\big) \\ 
=-\derlie_{\xt}(i\;w)+\jt\big(\derlie_{\xt}w\big) 
=-\derlie_{\xt}\big(\jt w\big)+\jt\big(\derlie_{\xt}w\big)=-\big(\derlie_{\xt}\jt\big)(w),
\end{multline*} 
for all $w\in\sections{\lt}$.
Thus,  
\[ \jtd v+\big(\derlie_{\xt}\jt\big)(v)=0, \qquad\forall v\in\sections{\lt} .\]
Since $\xt$ is a real section of $E$ and the endomorphisms $\jt$ and $\jtd$ are real, we also have 
\[ \jtd \cc{v}+\big(\derlie_{\xt}\jt\big)(\cc{v})=0, \qquad\forall \cc{v}\in\sections{\cc{\lt}}.\]
Therefore, $\jtd+\derlie_{\xt}\jt=0$, for all $t$. 
\end{proof}

\begin{thm}
\label{thm:main}
Assume that $M$ is a compact manifold, and $E$ is a 
Courant algebroid over $M$. Assume that $J_0$ and
$J_1$ are two generalized complex structures on $E$, which
are connected by a smooth family of generalized complex structures
 $t\mapsto\jt$ ($t\in[0,1]$) on $E$.
Assume there exists a smooth path $t\mapsto\bt$ in $\sections{\lt^*}$ 
such that $\otd=d_{\lt}\bt$, $\forall t$.
Then $J_0$ and
$J_1$ are isomorphic. 
\end{thm}


\begin{proof} 
The idea of the proof is to construct a smooth one-parameter 
family of automorphisms $t\mapsto\phi_t$ in $\Aut(E)$ such that 
$\phi_t\rond J_0=\jt\rond\phi_t$, $\forall t$.
Since $\lt^*\isomorphism\ltb$, there exists a unique time-dependent section $\zt\in\sections{\lt}$ such that 
$\bt=\ip{\ztb}{\cdot}$. Set $\xt=\tfrac{1}{2i}(\zt-\ztb)\in\sections{E}$. 
It follows from Lemma \ref{lem:bt} that \[ \jtd+\derlie_{\xt}\jt=0, \qquad \forall t .\] 
Therefore \[ \left.\tfrac{d}{d\tau} \bphi_{\tau}\inv\rond J_{\tau}\rond\bphi_{\tau}
\right|_{\tau=t}=\bphi_t\inv\rond\big(\jtd+\derlie_{\xt}\jt\big)\rond\bphi_t=0 ,\]
where $\bphi_t\in\Aut(E)$ is the flow of $\der_{\xt}\in\aut(E)$ (as in Section \ref{liederiv}), which always exists since $M$ is compact.
\end{proof}

\begin{rmk} When $M$ is not compact, the above theorem is still valid if $\xt$ is complete, i.e. if its flow exists for arbitrary time.
\end{rmk}

\section{Examples}

\subsection*{Example 1} \label{ex1}
A symplectic form on a Lie algebroid $A\to M$ is a smooth section $\Omega$ of $\wedge^2 A^*\to M$, which is $d_A$-closed and non-degenerate --- $a=0$ if $\Omega(a,b)=0$ for all $b\in\sections{A}$ \cite{\nest}. 
Here $d_A:\sections{\wedge^k A^*}\to\sections{\wedge^{k+1} A^*}$ is the Lie algebroid cohomology differential. 

Let $\Omega_0$, $\Omega_1$ be a pair of symplectic forms on $A$, which can be connected by a smooth family $t\mapsto\Omega_t$ of symplectic 2-forms. 
Having endowed $A^*$ with the trivial Lie algebroid structure, the pair $(A,A^*)$ becomes a Lie bialgebroid. 
On its "double", the Courant algebroid $E=A\oplus A^*$, we get the one-parameter family of generalized complex structures 
\[ J_t=\begin{pmatrix} 0 & -\pi\diese_t \\ \Omega\bemol_t & 0 \end{pmatrix} .\]
Here $\pi\diese_t:A^*\to A$ denotes the inverse of the bundle map $\Omega\bemol_t:A\to A^*$ associated to the symplectic form $\Omega_t$. 

It is well known that 
\[ \lambda_t:A\cmplx\to\lt:X\mapsto X-i\Omega\bemol_t X \] 
is a Lie algebroid isomorphism. 
Therefore, the diagram 
\[ \xymatrix{ \sections{\wedge^k \lt^*} \ar[d]_{d_{\lt}} \ar[r]^{\lambda_t^*} 
& \sections{\wedge^k A^*} \ar[d]^{d_A} \\ 
\sections{\wedge^{k+1} \lt^*} \ar[r]_{\lambda_t^*} & \sections{\wedge^{k+1} A^*} } \] 
is commutative.

Observe that $\lambda_t^*\otd=-\Otd$.
Indeed, 
\[ \lambda_t^*\otd(X,Y)=\ip{\lambda_t X}{\jtd\lambda_t Y}=
\left< \begin{pmatrix} X \\ -i\Ot\bemol X \end{pmatrix} , 
\begin{pmatrix} 0 & -\ptd\diese \\ \Otd\bemol & 0 \end{pmatrix} 
\begin{pmatrix} Y \\ -i\Ot\bemol Y \end{pmatrix} \right> 
= -\Otd(X,Y) .\]
Here we have used the identity $\Otd\bemol=-\Ot\bemol\rond\ptd\diese\rond\Ot\bemol$, which is easily obtained from 
$\pi_t\diese\rond\Ot\bemol=\idn$. 
Hence \[ d_A\Otd=-d_A\lambda_t^*\otd=-\lambda_t^* d_{\lt}\otd=0 .\]

Now assume there exists a family of smooth sections $\xi_t\in\sections{A^*}$ such that $\Otd=-d_A\xi_t$. 
It is easy to check that $\xi_t=\lambda_t^*\beta_t$, where $\beta_t=\ip{\cc{i\lambda_t\pi_t\diese\xi_t}}{\cdot}$.
Hence 
\[ \lambda_t^*\otd=-\Otd=d_A\xi_t=d_A\lambda_t^*\beta_t=\lambda_t^* d_{\lt}\beta_t ,\] i.e. $\otd=d_{\lt}\beta_t$. 

By Theorem \ref{thm:main}, we have $\bphi_{t*}J_t=J_0$, where $\bphi_t$ is the Courant algebroid automorphism generated by $\der_{\pi_t\diese\xi_t}$. Since $\pi_t\diese\xi_t\in\sections{A}$, $\bphi_t$ is actually induced by a Lie algebroid automorphism of $A$, which will also be denoted by $\bphi_t$ by abuse of notations. 

We have proved the following Lie algebroid version of Moser's Lemma \cite{\moser}.

\begin{cor}
Assume that $M$ is a compact manifold, and $A$ is a
Lie algebroid over $M$. Assume that $\Omega_0$ and 
$\Omega_1 \in \sections{\wedge^2 A^*}$ are two symplectic
structures on the Lie algebroid $A$
 which can be connected by a smooth family of  symplectic 
 structures
 $t\mapsto\Omega_t$ ($t\in[0,1]$) on $A$.
Assume there exists a smooth path $t\mapsto\xi_t$ in $\sections{A^*}$
such that $\Otd=-d_{A}\xi_t$, $\forall t$.
Then $\Omega_0$ and
$\Omega_1$ are isomorphic. That is, there is a Lie algebroid
automorphism $\phi$, isotopic to the identity, such that $\phi^*\Omega_1=\Omega_0$.
\end{cor}

\subsection*{Example 2} \label{ex2}

Let $A\to M$ be a Lie algebroid endowed with an integrable almost complex structure, i.e. a bundle map 
\[ \xymatrix{ A \ar[r]^j \ar[d] & A \ar[d] \\ M \ar[r]_{\idn} & M } \]
whose Nijenhuis tensor vanishes and such that $j^2=-\idn$. 
Its complexification $A\cmplx=A\otimes\CC$ decomposes as the direct sum $\aoz\oplus\azo$ of the eigenbundles of $j$ with eigenvalues $+i$ and $-i$ respectively.
Both $\aoz$ and $\azo$ are Lie subalgebroids of $A\cmplx$. 
The vector bundle $\aoz$ is naturally an $\azo$-module: the representation of $\azo$ on $\aoz$ is the map 
\[ \sections{\azo}\times\sections{\aoz}\to\sections{\aoz}:(X,Y)\mapsto X\action Y:=\proz\lie{X}{Y} ,\]
where $\proz$ stands for the canonical projection $A\cmplx=\azo\oplus\aoz\to\aoz$. 

Now consider the dual vector bundle $A^*\to M$ as a trivial Lie algebroid --- the anchor and the bracket are both zero. 
The pair $(A,A^*)$ is a Lie bialgebroid and the endomorphism 
\[ J = \begin{pmatrix} -j & 0 \\ 0 & j^* \end{pmatrix} \] 
is a generalized complex structure on its double $E:=A\oplus A^*$.
The $+i$-eigenbundle $L$ of $J$ is naturally isomorphic, as a Lie algebroid, to the semi-direct product Lie algebroid $\azo\ltimes(\aoz)^*$, where $(\aoz)^*$ is the $\azo$-module dual to $\aoz$. 

The following lemma can be easily verified.

\begin{lem} \label{lem:shift}
Let $A$ be a Lie algebroid, $M$ be an $A$-module and $A\ltimes M^*$ be the resulting semi-direct product Lie algebroid. Then the diagram 
\[ \xymatrix{ \cochains{k-1}{A}{M} \ar[r]^{\alpha} \ar[d]_{D} & \cochains{k}{A\ltimes M^*}{\CC} \ar[d]^{D} \\ 
\cochains{k}{A}{M} \ar[r]_{\alpha} & \cochains{k+1}{A\ltimes M^*}{\CC} } \]
commutes for all $k\ge 0$.
Here $\cochains{k}{A}{M}:=\sections{\wedge^k A^*\otimes M}$, the $D$'s are the coboundary operators for the Lie algebroid cohomologies of $A$ (with values in the module $M$) and $A\ltimes M^*$, while $\alpha$ is the skew-symmetrization map 
\[ (\wedge^k A^*)\otimes M \to (\wedge^k A^*)\wedge M \inclusion \wedge^{k+1}(A\ltimes M^*)^* .\]
\end{lem}

Let $t\mapsto j_t,\;t\in[0,1]$ be a smooth family of integrable almost complex structues on $A$. 
And let \[ \jt=\begin{pmatrix} -j_t & 0 \\ 0 & j_t^*\end{pmatrix} \] be the associated family of generalized complex structures on the Courant algebroid $E=A\oplus A^*$. 

Taking $A=\azot$ and $M=\aozt$ in Lemma~\ref{lem:shift}, we get that 
\[ \xymatrix{ \cochains{0}{\azot}{\aozt} \ar[r]^{\alpha} \ar[d]_{\delbar_t} 
& \cochains{1}{\lt}{\CC} \ar[d]^{d_{L_t}} \\ 
\cochains{1}{\azot}{\aozt} \ar[r]_{\alpha} & \cochains{2}{\lt}{\CC} } \]
is commutative for all $t\in [0,1]$. 
Here $\delbar$ denotes the differential operator of the Lie algebroid cohomology of $\azot$ with values in the module $\aozt$. 

The time-derivative $\dot{j}_t$ may be seen as an element of 
$\sections{(\azot)^*\otimes\aozt}=\cochains{1}{\azot}{\aozt}$ because $\dot{j}_t$ swaps $\azot$ and $\aozt$. 
Moreover, $\otd=-\alpha(\dot{j}_t)$. Indeed, for all $X,Y\in\sections{\azot}$ and all $\xi,\eta\in\sections{(\aozt)^*}$, one has 
\begin{multline*} \otd(X+\xi,Y+\eta)=\ip{X+\xi}{\jtd(Y+\eta)}=\ip{X+\xi}{\dot{j}_t Y-\dot{j}_t^*\eta} \\ =\thalf\big(\xi(\dot{j}_t Y)-\eta(\dot{j}_t X)\big)=-\alpha(\dot{j}_t)(X+\xi,Y+\eta) .\end{multline*}

Assume there exists a family $t\in [0,1]\mapsto z_t\in\sections{\azot}$ such that \[ \delbar_t(\thalf\cc{z_t})=-\dot{j}_t .\] (Remark that $\thalf\cc{z_t}\in\sections{\aozt}=\cochains{0}{\azot}{\aozt}$.) 
Then \[ \otd=-\alpha(\dot{j}_t)=\alpha\rond\delbar_t(\thalf\cc{z_t})
=\dlt\rond\alpha(\thalf\cc{z_t})=\dlt(\ip{\cc{z_t}}{\cdot}) ,\]
since $\ip{\cc{z_t}}{X+\xi}=\thalf\xi(\cc{z_t})=\thalf\alpha(\cc{z_t})\big(\xi\big)$ for all $X+\xi\in\lt=\azot\ltimes(\aozt)^*$. 

Set $x_t=\tfrac{1}{2i}(z_t-\cc{z_t})\in\sections{A}$. According to Theorem \ref{thm:main}, we have $\bphi_{t*}J_t=J_0$, where $\bphi_t$ is the flow of automorphisms of the Courant algebroid $E=A\oplus A^*$ generated by $\der_{x_t}$. 
Here, since the section $x_t$ of $E$ actually lies in the Lie algebroid $A$, $\bphi_t$ is induced by a family of Lie algebroid isomorphisms. 

We have proved 

\begin{cor} 
Let $A$ be a Lie algebroid over a compact manifold $M$. Let $j_0$ and $j_1$ be two integrable almost complex structures on $A$, which are connected by a smooth family $j_t$ of integrable almost complex structures. Assume there exists a smooth family $z_t\in\sections{\azot}$ such that $\delbar_t \cc{z_t}=-2\dot{j}_t$. Then $j_0$ and $j_1$ are isomorphic, i.e. there exists a Lie algebroid automorphism $\bphi$ such that $j_1=\bphi_* j_0 (:=\bphi\inv\rond j_0\rond\bphi)$. 
\end{cor}

When $A$ is the tangent Lie algebroid $TM$, this is a 
classical result \cite{\kodaira}. 

\section{Hamiltonian operators}

Let $J$ be a generalized complex structure on a Courant algebroid $E$. Then $E=L\oplus\cc{L}$ and $L^*\isomorphism\cc{L}$, where $L$ is the $+i$-eigenbundle of $J$. And $(L,L^*)$ is a complex Lie bialgebroid. 

Modulo the isomorphism $E\to E^*:e\mapsto\ip{e}{\cdot}$, which identifies $\cc{L}$ with $L^*$ (resp. $L$ with $\cc{L}^*$), any section $H\in\sections{\wedge^2L^*}$ (resp. $\cc{H}\in\sections{\wedge^2\cc{L}^*}$) can be seen as a map $H':L\to\cc{L}$ (resp. $\cc{H}':\cc{L}\to L$). 
The graph $L_H=\gendex{(v,H'v)}{v\in L}\subset L\oplus\cc{L}$ of $H':L\to\cc{L}$ is a Dirac structure of $E$ if, and only if, the Maurer-Cartan equation 
\[ \dl H+\thalf\lie{H}{H}_{L^*}=0 \]
holds \cite{\liu}.
In this case, $H$ is called a \emph{Hamiltonian operator}. 
Here $\dl:\sections{\wedge^{\bullet} L^*}\to\sections{\wedge^{\bullet+1} L^*}$ is the Lie algebroid differential of $L$ and $\lie{\cdot}{\cdot}_{L^*}:\sections{\wedge^k L^*}\otimes\sections{\wedge^l L^*}\to\sections{\wedge^{k+l-1} L^*}$ is the Schouten bracket of the Lie algebroid $L^*$. 

When $H$ is a Hamiltonian operator, $L_H$ is a Lie algebroid. 
The vector bundle isomorphism 
\[ L\to L_H:v\mapsto(v,H'v) \] 
induces the cochain complex isomorphism 
\beq{eq:cci} 
\xymatrix{ \dots \ar[r]^{d_{L_H}} & \sections{\wedge^{\bullet} L_H^*} \ar[r]^{d_{L_H}} \ar[d]^{\isomorphism} & 
\sections{\wedge^{\bullet+1} L_H^*} \ar[r]^{d_{L_H}} \ar[d]^{\isomorphism} & \dots \\ 
\dots \ar[r]^{d_H} & \sections{\wedge^{\bullet} L^*} \ar[r]^{d_H} & 
\sections{\wedge^{\bullet+1} L^*} \ar[r]^{d_H} & \dots }
\eeq
where $d_H=\dl+\lie{H}{\cdot}_{L^*}$ and $d_{L_H}$ is the Lie algebroid coboundary operator of $L_H$. 

\begin{lem}\label{lem:invertible}
Assume $H\in\sections{\wedge^2 L^*}$ is a Hamiltonian operator. 
Then $L_H$ is the $+i$-eigenbundle of a generalized complex structure on $E$ if, and only if, the bundle map $\cc{H}'\rond H'-1:L\to L$ is invertible. 
\end{lem} 

\begin{proof}
First, note that $\cc{(L_H)}=(\cc{L})_{\cc{H}}$. 
Hence $(v,H'v)\in L_H\cap\cc{L_H}$ iff $\exists w\in\cc{L}$ such that $(v,H'v)=(\cc{H}'w,w)$, which holds if, and only if, $(\cc{H}'\rond H'-1)v=0$. It thus follows that $L_H\cap\cc{L_H}=\{0\}$ iff $\cc{H}'\rond H'-1$ is invertible.
\end{proof}

\begin{ex}
Consider the generalized complex structure $J=\left( \begin{smallmatrix} -j & 0 \\ 0 & j^*
 \end{smallmatrix} \right)$ on the standard
Courant algebroid $E=\TTM$ corresponding to a
complex structure $j: TM\to TM$. Then we have $L=T^{0, 1}\ltimes (T^{1, 0})^*$ and $\bar{L}=T^{1, 0}\ltimes (T^{0, 1})^*$, where the Lie algebroid
structure on $L$ is the semi-direct product, similarly for $\bar{L}$.
Let $\pi\in\sections{\wedge^2 T^{1, 0}}$. It is simple to see that
$\cc{\pi}' \rond \pi' =0$. Hence $\cc{\pi}' \rond \pi' -
1=-1$ is invertible. On the other hand,  $\pi$ satisfies 
the Maurer-Cartan equation if and only if $\delbar\pi=0$ and
$\schouten{\pi}{\pi}=0$. That is, $\pi$ is a holomorphic Poisson.
Thus we recover the well known fact that a holomorphic Poisson
is a generalized complex structure \cite{\gualtieri,\crainic} on the standard Courant
algebroid $\TTM$.
\end{ex}

Now, let $t\in[0,1]\mapsto\Ht\in\sections{\wedge^2 L^*}$ be a family of Hamiltonian operators such that 
$\cc{\Ht}'\rond\Ht'-1$ is invertible for all $t\in[0,1]$. 
By Lemma \ref{lem:invertible}, we obtain a family of generalized complex structures $t\mapsto\jt$. 
In this case, its corresponding $\otd$ has a very simple description. 

\begin{lem}\label{lem:otdhtd}
$\Theta_t^*\otd=-i\Htd \quad \in\sections{\wedge^2 L^*}$
\end{lem} 

Here $\Theta_t$ is the vector bundle map 
\[ L\to E=L\oplus\cc{L}:v\mapsto(v,\Ht'v) \] 
identifying $L$ with $L_{\Ht}$. 

\begin{proof}
By definition, we have 
\[ \jt(v,\Ht'v)=i(v,\Ht'v) .\]
Differentiating w.r.t. $t$, we get 
\[ \jtd(v,\Ht'v)+\jt(0,\Htd'v)=i(0,\Htd'v) .\]
It thus follows that 
\begin{multline*} 
\Theta_t^*\otd(v,w)=\ip{(v,\Ht'v)}{\jtd(w,\Ht'w)}=\ip{(v,\Ht'v)}{(i-\jt)(w,\Ht'w)} \\ 
=\ip{(i+\jt)(v,\Ht'v)}{(w,\Ht'w)}=2i\ip{(v,\Ht'v)}{(w,\Ht'w)}=-i\Htd(v,w) 
.\end{multline*}
\end{proof}

From \eqref{eq:cci} and Lemma \ref{lem:otdhtd}, it follows that \eqref{eq:wt} is equivalent to 
\[ \frac{d\Ht}{dt}=\dl\yt+\lie{\Ht}{\yt}_{L^*} ,\] where $\yt=i\Theta_t^*\beta_t\in\sections{L^*}$. 
We have proved

\begin{thm} 
Let $E$ be a Courant algebroid, over a compact manifold $M$, whose complexified admits a decomposition $E\cmplx=L\oplus\cc{L}$ as the direct sum of a pair of complex conjugate Lie algebroids. And let $\Ht\in\sections{\wedge^2 L^*}$ be a smooth family of Hamiltonian operators such that 
$\cc{\Ht}'\rond\Ht'-1$ is invertible. 
Assume there exists a smooth family $\yt\in\sections{L^*}$ solving the equation 
\beq{eq:dHt} \frac{d\Ht}{dt}=\dl\yt+\lie{\Ht}{\yt}_{L^*} .\eeq
Then the generalized complex structures associated to $H_0$ and $H_1$ are isomorphic. 
\end{thm}

\begin{ex}
Take $L=T^{0,1}\ltimes(T^{1,0})^*$ as in Section \ref{ex2}. 
Let $\pi_t\in\sections{\wedge^2 T^{1,0}}$ be a family of holomorphic 
Poisson bivector fields. It is simple to check that 
$\yt=X_t^{1,0}+\xi_t^{0,1}\in\cc{L}=T^{1,0}\ltimes(T^{0,1})^*$ satisfies \eqref{eq:dHt} iff
\[ \delbar X_t^{1, 0}=0, \qquad \bar{\partial}\xi_t^{0,1}=0 \quad \text{and} \quad 
\dot{\pi}_t=\lie{\pi_t}{X_t^{1,0}} .\]
In particular, one may take $\xi_t^{0,1}=0$.
Thus we obtain the following standard result:
if there exists a family of complete holomorphic vector fields
$X_t^{1, 0}$ satisfying the equation
\[ \dot{\pi}_t=\lie{\pi_t}{X_t^{1, 0}} ,\]
then the holomorphic Poisson structures $\pi_0$ and $\pi_1$ are isomorphic.
Note that, in this case, $x_t=2\RE X_t^{1,0}\in\vf(M)$ and $x_t$ is complete
iff $X_t^{1,0}$ is complete.
\end{ex}

\section{Darboux-Weinstein theorem}

As an application of Theorem \ref{thm:main}, we present
a Darboux-Weinstein style theorem for generalized complex structures.

Let $J$ be a generalized complex structure on the standard Courant algebroid $\TTM\to M$. 
Then $\sections{\TTM}$ is a module over the ring $\cinf(M)$. Given a point $o\in M$, choose local coordinates $(x_1,\dots,x_n)$ identifying an open neighbourhood $U$ of $o$ in $M$ with a unit ball in $\RR^n$ centered at the origin (and mapping $o$ to $0$). Then $\gendex{e_1=\partial_{x^1},\dots,e_n=\partial_{x^n},e_{n+1}=dx^1,\dots,e_{2n}=dx^n}{}$ is a local frame of $\TTM$ over $U$, which enjoys the following remarkable properties: 
\begin{gather} 
\courant{e_i}{e_j}=0, \qquad \forall i,j\in\gendex{1,\dots,2n}{}; \label{eq:RP1} \\ 
\ip{e_i}{e_j} \text{ is constant on } U, \qquad \forall i,j\in\gendex{1,\dots,2n}{}. \label{eq:RP2} 
\end{gather} 

Now consider the matrix representation of $J$ in this local frame: 
\[ J \ejx=\sum_i a^i_j(x) \eix, \qquad x\in U. \]
The $a^i_j$'s are smooth functions on $U\subset\RR^n$.
Set 
\beq{eq:J_t} \jt \ejx=\sum_i a^i_j(tx) \eix, \eeq 
for all $t\in [0,1]$ and $x\in U$.
Clearly, $\jt^2=-\idn$, $J_1=J$ and $J_0$ is represented by the constant matrix $a^i_j(0)$. 
Moreover, \eqref{eq:RP1} implies that $J_0$ is integrable. 

\begin{lem} 
All $\jt$'s are generalized complex structures on the restriction of $\TTM$ to $U$. 
\end{lem}

\begin{proof}
Let $\gendex{e^k}{k=1,\dots,2n}$ be the dual frame to $\gendex{e_i}{i=1,\dots,2n}$. 
From a tedious computation using the relations \eqref{eq:RP1}, \eqref{eq:RP2} and \[ \courant{fA}{gB}=fg\courant{A}{B}+f\big(\anchor(A)g\big)B-g\big(\anchor(B)f\big)A +\ip{A}{B}(g\DD f-f\DD g) \] 
(where $f,g\in\cinf(M)$ and $A,B\in\sections{\TTM}$), it follows that 
\begin{gather*} 
\big(e^k,\courant{\jt e_i}{\jt e_j}\big)(x)=t\big(e^k,\courant{J e_i}{J e_j}\big)(tx), \\ 
\big(e^k,\courant{\jt e_i}{e_j}\big)(x)=t\big(e^k,\courant{J e_i}{e_j}\big)(tx), \\ 
\intertext{and}
\big(e^k,\jt\courant{\jt e_i}{e_j}\big)(x)=t\big(e^k,J\courant{J e_i}{e_j}\big)(tx). 
\end{gather*}
Hence the integrability of $J$ at the point $tx$ implies the integrability of $\jt$ at the point $x$. Since $J$ is integrable on $U$, $U$ is the unit ball around $0$ in $\RR^n$ and $t\in [0,1]$, the conclusion follows.
\end{proof}

As an immediate consequence of Theorem \ref{thm:main}, we have 

\begin{cor} 
Assume there exists a smooth path $t\in [0,1]\mapsto \beta_t\in\sections{\lt^*|_U}$ such that 
$\otd=d_{\lt}\beta_t$. Then the restriction of $J$ to $U$ is isomorphic to the constant generalized complex structure $J_0$. 
\end{cor}

\begin{thm} \label{thm:7.6}
Let $J$ be a generalized complex structure (with $+i$-eigenbundle $L$) on the standard Courant algebroid $\TTM\to M$. Let $(x_1,\dots,x_n)$ be local coordinates on $M$ identifying an open domain $U$ with the unit ball centered at the origin in $\RR^n$. And define $\jt$ as in \eqref{eq:J_t}. 
If there exists $\beta\in\sections{L^*|_U}$ such that 
\[ d_L\beta=\dot{\omega}_1, \qquad \beta|_0=0 \quad \text{and} \quad d\beta|_0=0, \]
then $J|_U$ is isomorphic to a constant generalized complex structure. 
\end{thm}

\begin{rmk}
Using the Darboux and Newlander-Nirenberg theorems, a result of Gualtieri \cite[Theorem~4.35]{\gualtieri} can be restated as follows: 
if $m\in M$ is a regular point of $J$ (i.e. the corresponding Poisson tensor is
regular at $m$), then, locally around $m$, $J$ must be isomorphic to a constant generalized complex structure. Therefore the condition of Theorem~\ref{thm:7.6} must fail in a neighbourhood of a regular point.
\end{rmk}

The proof is splitted into several lemmas. 

\begin{lem} \label{lem:7.7}
Let $E_1$ and $E_2$ be two Courant algebroids over the same base manifold $M$, with anchor maps $\anchor_1$ and $\anchor_2$ respectively. Let 
\[ \xymatrix{ E_1 \ar[d] \ar[r]^{\bphi} & E_2 \ar[d] \\ M \ar[r]_{\vphi} & M } \] 
be an isomorphism of vector bundles such that 
\[ \anchor_2\rond\bphi=\vphi_*\rond\anchor_1 \] 
and \[ \courant{\bphi\inv\rond\alpha\rond\vphi}{\bphi\inv\rond\beta\rond\vphi}=\bphi\inv\rond\courant{\alpha}{\beta}\rond\vphi, 
\qquad \forall \alpha,\beta\in\sections{E_2} . \]
Then $\bphi$ establishes a 1-1 correspondence between the subbundles of $E_1$ and  $E_2$ which are Lie algebroids. 
\end{lem}

\begin{lem} 
The vector bundle map 
\[ \xymatrix{ \TTM|_U \ar[d] \ar[r]^{\bphi_t} & \TTM|_U \ar[d] \\ 
U \ar[r]_{\vphi_t} & U, } \]
where \[ \vphi_t(x_1,\dots,x_n)=(t x_1,\dots,t x_n) \] and 
\[ \bphi_t(\eix)=t\eitx ,\]
satisfies the assumptions of Lemma \ref{lem:7.7}.
Moreover, we have \[ J\rond\bphi_t=\bphi_t\rond\jt .\] 
\end{lem} 

\begin{cor} 
The "multiplication by $t$" map 
\[ \xymatrix{ L_t \ar[d] \ar[r]^{\bphi_t} & L \ar[d] \\ 
U \ar[r]_{\vphi_t} & U } \] 
is a Lie algebroid isomorphism. 
\end{cor} 

Therefore we have \[ \bphi_t^*\rond d_L=d_{\lt}\rond \bphi_t^* .\] 

\begin{lem} \label{lem:6.11}
\[ t^3\otd=\bphi_t^*\dot{\omega}_1 \]
\end{lem}

\begin{proof}
Differentiating \eqref{eq:J_t} w.r.t. $t$, we get 
\[ \jtd\ejx=\sum_{k,l}x^l\left.\frac{\partial a^k_j}{\partial x^l}\right|_{tx}\ekx .\]
Therefore 
\[ \otd(\eix,\ejx)=\ip{\eix}{\jtd\ejx}=b_{ij}x^l\left.\frac{\partial a^k_j}{\partial x^l}\right|_{tx} .\]
Substituting $1$ to $t$ and $tx$ to $x$ in the last relation, one gets 
\[ \dot{\omega}_1(\eitx,\ejtx)=b_{ij}tx^l\left.\frac{\partial a^k_j}{\partial x^l}\right|_{tx} .\] 
Now, multiplying by $t^2$, we obtain 
\[ \dot{\omega}_1(t\eitx,t\ejtx)=b_{ij}t^3x^l\left.\frac{\partial a^k_j}{\partial x^l}\right|_{tx} \] 
or, equivalently, 
\[ \dot{\omega}_1(\bphi_t\eix,\bphi_t\ejx)=t^3\otd(\eix,\ejx) .\]
\end{proof}

\begin{proof}[Proof of Theorem \ref{thm:7.6}]
According to Lemma \ref{lem:6.11}, 
\[ t^3\otd=\bphi_t^*\dot{\omega}_1=\bphi_t^*d_L\beta=d_{\lt}\bphi_t^*\beta .\]
Now set $\beta_t=t^{-3}\bphi_t^*\beta$. 
We have $\beta_t(\eix)=t^{-3}\beta(\bphi_t\eix)=t^{-2}\beta(\eitx)$, which shows that $\beta_t$ is well-defined for $t=0$ since $\beta|_0=0$ and $d\beta|_0=0$. 
The conclusion follows from Theorem \ref{thm:main} since $\otd=d_{\lt}\beta_t$.
\end{proof}

\begin{rmk}
When $J$ is the generalized complex structure corresponding to a symplectic
structure, it is straightforward to check that the Poincar\'e lemma for de Rham cohomology ensures that the conditions of Theorem~\ref{thm:7.6}
hold at any point $m\in M$. 
As a consequence, we recover the Darboux-Weinstein theorem.
However the one-parameter family of generalized complex structures (or family of
symplectic structures) considered here is different from the one used in \cite{\weinstein}.

On the other hand, when $J$ is the generalized complex structure 
corresponding to an integrable complex structure
$j:TM\to TM$, Kohn gave a proof of the Newlander-Nirenberg theorem using
exactly the same $j_t$ as ours \cite{\kohn}.
Kohn's essential idea and major achievement was to prove the Poincar\'e lemma for Dolbeault cohomology using the description of this cohomology in terms of the integrable almost complex structure rather than holomorphic charts. 
As a consequence, he was able to deform the local chart used to define the family of local almost complex structures $t\mapsto j_t$ joining the given $j$ to the constant $j_0$ into a local chart in which $j$ is constant. 
It would be interesting to compare our condition with the one in \cite{\kohn}
in that case.
\end{rmk}

\begin{rmk} 
Note that the degree 2 Lie algebroid cohomology of $L$
may not always be locally trivial for generic generalized complex structures 
even though the Poincar\'e Lemma holds for
both the symplectic and complex cases. For instance, 
when $J$ corresonds to a holomorphic Poisson struture $\pi$,
the Lie algebroid cohomology of $L$ is simply the Poisson
cohomology of $\pi$ \cite{\hp}, which
may not always be locally trivial.

On the other hand, the condition of Theorem~\ref{thm:7.6} 
is stronger than just requiring $H^2(L|_U,\CC)$ 
to vanish locally on a neighbourhood $U$.
For instance, when $J$ is the generalized complex structure 
corresponding to the holomorphic Lie Poisson structure of a complex
simple Lie algebra, $H^2(L|_U,\CC)$ vanishes on a nbd $U$
of $0$. However $J$ cannot be locally trival on neighbourhoods of $0$.
\end{rmk}


\bibliographystyle{hamsplain}
\bibliography{moserbiblio}

\end{document}